\documentclass[12pt]{article}
\usepackage{mathrsfs}

\usepackage{t1enc}
\usepackage[latin1]{inputenc}
\usepackage[english]{babel}
\usepackage{amsmath,amsthm,amssymb}
\usepackage{amsfonts}
\usepackage{latexsym}
\usepackage{enumerate}

\usepackage{graphicx}
\usepackage[natural]{xcolor}
\usepackage{tikz}
\usepackage{tikz-3dplot}
\usetikzlibrary{shapes}
\usetikzlibrary{arrows,decorations.pathmorphing,backgrounds,positioning,fit,petri,automata}
\usetikzlibrary {positioning}

\theoremstyle{plain}
\newtheorem{thm}{Theorem}

\newtheorem{prob}[thm]{Problem}
\newtheorem{coro}[thm]{Corollary}

\newtheorem{claim}{Claim}

\theoremstyle{plain}

\theoremstyle{plain}

\theoremstyle{plain}

\title{Nonempty intersection of longest paths in graphs without forbidden pairs}
\author{Yuping Gao$^{a}$, Songling Shan$^{b}$\\
{\small a. School of Mathematics and Statistics, Lanzhou University, Lanzhou 730000, China}\\
{\small b. Department of Mathematics, Illinois State University, Normal, IL 61790, USA}}
\date{}

\begin{document}
\baselineskip 0.65cm

\maketitle
\begin{abstract}

In 1966, Gallai asked whether all longest paths in a connected graph have a nonempty intersection. The answer to this question is not true in general  and various counterexamples have been found. However, there is a  positive solution to Gallai's question for many well-known classes of graphs
such as split graphs, series parallel graphs, and $2K_2$-free graphs.  Among all the graph classes  that support Gallai's question, almost all of them were shown to be Hamiltonian under certain conditions. This observation motivates  us to investigate Gallai's question in graphs that are ``close'' to Hamiltonicity properties. Let  $\{R,S\}$ be a pair of  connected graphs.  In particular, in this paper, we show that Gallai's question  is affirmative for all connected $\{R,S\}$-free graphs such that every 2-connected  $\{R,S\}$-free graph is Hamiltonian. These pairs $\{R,S\}$ were completely  characterized in 1990s.

\medskip

\end{abstract}

\section{Introduction}
All graphs considered in this paper are connected, undirected and simple. A path in a graph is a \emph{longest path} if there exist no other paths in the graph that are strictly longer. It is well known that any two longest paths share a common vertex in any connected graph. In 1966, Gallai asked whether all longest paths in a connected graph have a common vertex \cite{Gallai1966}. The answer to this question is known to be negative and the first counterexample was given by Walther on a graph with 25 vertices in \cite{Walther1969}. The smallest graph answering Gallai's question negatively is a graph on 12 vertices, found by Walther and Voss in \cite{WV1974} and independently by Zamfirescu in \cite{Zamfirescu1976} (see Fig. 1). Brinkmann and Van Cleemput~\cite{BV2012} proved that there is no counterexample to Gallai's question with  less than 12 vertices.

\begin{center}
\begin{tikzpicture}

{\tikzstyle{every node}=[draw,circle,fill=black,minimum size=4pt,
                            inner sep=0pt]
 \draw (1,3) node (v1)  {}
        -- ++(300:2cm) node (v2)  {}
        -- ++(180:2cm) node (v3)   {}
        --(v1) ;
 \draw (1.6,2) node (v4)  {}
        -- ++(300:2cm) node (v5)  {}
        -- ++(180:2cm) node (v6)   {}
        --(v4) ;
 \draw (0.4,2) node (v7)  {}
        -- ++(300:2cm) node (v8)  {}
        -- ++(180:2cm) node (v9)   {}
        --(v7) ;
 \draw (1,4) node (v10)  {}
        --(v1) ;
 \draw (3.3,-0.5) node (v11)  {}
        --(v5) ;
 \draw (-1.2,-0.5) node (v11)  {}
        --(v9) ;
   }

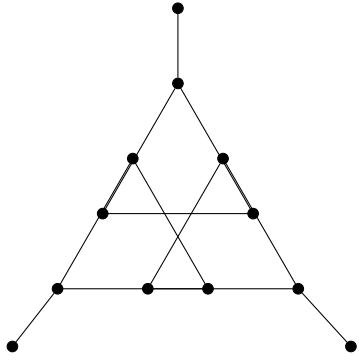
\node at (1,-1.5) {Fig. 1: Counterexample of Walther, Voss and Zamfirescu};
\end{tikzpicture}
\end{center}

Gallai's question is true for many classes of graphs, such as split graphs \cite{KP1990}, series-parallel graphs \cite{CEFHSYY2017}, graphs with matching number at most 3 \cite{Chen2015}, $2K_2$-free graphs \cite{GS2018}, etc.  Observe that all these graphs are either shown to be Hamiltonian under certain conditions or have ``strong'' Hamiltonicity properties, i.e., they have many long paths and long cycles. This observation motivates us to investigate Gallai's question in graphs that are ``close'' to have a Hamiltonian path or Hamiltonian cycle.

For a connected graph $H$, a graph $G$ is said to be \emph{$H$-free} if $G$ does not contain $H$ as an induced subgraph. For a set of connected graphs $\mathcal{H}$, $G$ is said to be \emph{$\mathcal{H}$-free} if $G$ is $H$-free for every $H\in \mathcal{H}$. Specifically, call $\mathcal{H}$ a \emph{forbidden pair} if $|\mathcal{H}|=2$, and if $\mathcal{H}=\{R,S\}$, we simply
write that $G$ is $(R,S)$-free.

 Let $P_{n}$ denote a  path of order $n$.
For nonnegative integers $k$, $\ell$ and $m$, let $N_{k,\ell,m}$ be a graph obtained from $K_3$ and three vertex-disjoint paths $P_{k+1}$, $P_{\ell+1}$, $P_{m+1}$ by identifying each of the vertices of the $K_{3}$ with one endvertex of one of the paths.  Let $Z_2=N_{2,0,0}, Z_3=N_{3,0,0},  B_{1,1}=N_{1,1,0}$, and
$B_{1,2}=N_{1,2,0}$.

Recently, Cerioli and Lima \cite{CERIOLI2019}  showed that Gallai's question is true for several classes of graphs including  chain graphs,  $P_4$-sparse graphs, starlike graphs, and $(K_{1,3}, P_5)$-free graphs.  In fact, $(K_{1,3}, P_5)$ is one pair of the graphs whose exclusion forces a 2-connected graph to be Hamiltonian. We call a graph {\it traceable\/} if it has a Hamiltonian path.
In 1997, Faudree  and Gould~\cite{FG1997}  proved the following result.

\begin{thm}[\cite{FG1997}]\label{traceable}  Let $R$ and $S$ be connected graphs $(R,S\neq P_3)$ and let $G$ be a connected graph. Then $G$ is $(R,S)$-free implies that $G$ is traceable if and only if $R=K_{1,3}$ and $S$ is one of the graphs $C_3,P_4,Z_1,B_{1,1}$ or $N_{1,1,1}$.
\end{thm}

For 2-connected graphs, Bedrossian in 1991~\cite{Bedrossian1991} obtained all the forbidden pairs for
Hamiltonian cycles as below.

\begin{thm}[\cite{Bedrossian1991}]\label{forbiddenpairH1}Let $R$ and $S$ be connected graphs $(R,S\neq P_3)$ and $G$ a $2$-connected graph. Then $G$ is $(R,S)$-free implies that $G$ is Hamiltonian if and only if $R=K_{1,3}$ and $S$ is one of the graphs $C_3,P_4,P_5,P_6,Z_1,Z_2,B_{1,1},B_{1,2}$ or $N_{1,1,1}$.
\end{thm}

The following result indicates that $(K_{1,3},Z_{3})$ is an additional forbidden pair for Hamiltonian cycles  if we consider  graphs with  large order.

\begin{thm}[\cite{FG1997}]\label{forbiddenpairH}  Let $R$ and $S$ be connected graphs $(R,S\neq P_3)$ and $G$ a $2$-connected graph of order $n\geq 10$. Then $G$ is $(R,S)$-free implies that $G$ is Hamiltonian if and only if $R=K_{1,3}$ and $S$ is one of the graphs $C_3,P_4,P_5,P_6,Z_1,Z_2,Z_3,B_{1,1},B_{1,2}$ or $N_{1,1,1}$ \emph{(}See Fig. \emph{2)}.
\end{thm}

\begin{center}
\begin{tikzpicture}

{\tikzstyle{every node}=[draw,circle,fill=black,minimum size=4pt,
                            inner sep=0pt]
      \node at (0,2) {};
  \foreach \numbera in {1,...,3}{
      \node[draw,circle]  (M-\numbera) at (2,\numbera) {};
    }
  \foreach \numbera in {1,...,3}{
    \path (0,2) edge (M-\numbera);
  }
}
\node at (1,0.5) {$K_{1,3}$};

{\tikzstyle{every node}=[draw,circle,fill=black,minimum size=4pt,
                            inner sep=0pt]
 \draw (4,3) node (v1)  {}
        -- ++(300:2cm) node (v2)  {}
        -- ++(180:2cm) node (v3)   {}
        -- (v1)
        ;
   }
\node at (4.1,0.5) {$C_{3}$};

{\tikzstyle{every node}=[draw,circle,fill=black,minimum size=4pt,
                            inner sep=0pt]
    \draw (6,2) node (v1)  {}
        -- ++(0:1cm) node (v2){}
        -- ++(0:1cm) node (v3){}
        -- ++(0:1cm) node (v4){}
        -- ++(0:1cm) node (v5){}
        -- ++(0:1cm) node (v6){}
        ; }
\node at (8.5,1) {$P_{6}$};
\end{tikzpicture}
\end{center}

\begin{center}
\begin{tikzpicture}

{\tikzstyle{every node}=[draw,circle,fill=black,minimum size=4pt,
                            inner sep=0pt]
 \draw (0,3) node (v1)  {}
        -- ++(180:2cm) node (v2)   {}
         -- ++(300:2cm) node (v3)  {}
         -- ++(0:1cm) node (v4){}
         -- ++(0:1cm) node (v5){}
         -- ++(0:1cm) node (v6){};
        \draw (v1)--(v3) ;
   }
\node at (0,0.5) {$Z_{3}$};

{\tikzstyle{every node}=[draw,circle,fill=black,minimum size=4pt,
                            inner sep=0pt]
 \draw (5,3) node (v1)  {}
        -- ++(180:2cm) node (v2)   {}
         -- ++(300:2cm) node (v3)  {}
         -- ++(270:1cm) node (v4)  {};
 \draw (v1)--(v3) ;

\draw (2,3) node (v5)  {}
         -- (v1);
\draw (6,3) node (v6)  {}
         -- (v2);
   }
\node at (4,-0.5) {$N_{1,1,1}$};

{\tikzstyle{every node}=[draw,circle,fill=black,minimum size=4pt,
                            inner sep=0pt]
 \draw (9,3) node (v1)  {}
        -- ++(180:2cm) node (v2)   {}
         -- ++(300:2cm) node (v3)  {}
         -- ++(0:1cm) node (v4){}
         -- ++(0:1cm) node (v5){};
        \draw (v1)--(v3) ;

 \draw (10,3) node (v6)  {}
         -- (v1);
   }
\node at (8.5,0.5) {$B_{1,2}$};
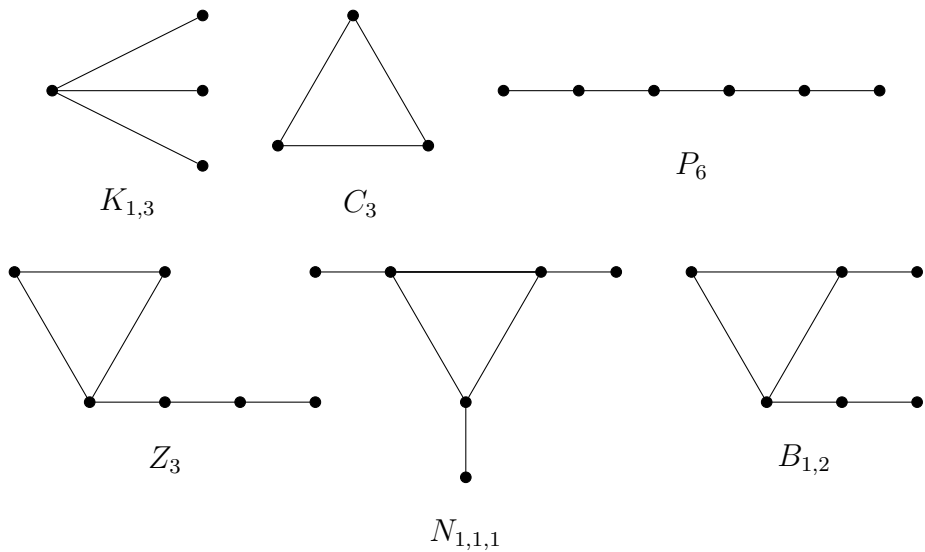
\node at (4,-1.5) {Fig. 2: Forbidden pairs for Hamiltonian cycles};
\end{tikzpicture}
\end{center}

By Theorem~\ref{forbiddenpairH1}, $(K_{1,3}, Z_3)$ is the only forbidden pair that requires $G$
to have order at least 10. In fact, Faudree  and Gould~\cite{FGRS1995} characterized all 2-connected
$(K_{1,3}, Z_3)$-free graphs that are not Hamiltonian, as listed below.

\begin{thm}[\cite{FGRS1995}] \label{CZ3}  If $G$ is a $2$-connected $(K_{1,3},Z_3)$-free graph, then $G$ is either Hamiltonian or isomorphic to $H_1$ or $H_2$ $($see Fig. $3)$.
\end{thm}

\begin{center}
	\begin{tikzpicture}
	{\tikzstyle{every node}=[draw,circle,fill=white,minimum size=4pt,
		inner sep=0pt]
		\draw (0,3) node (v1)  {$v_1$}
		-- ++(0:3cm) node (v2)  {$v_2$}
		-- ++(270:3cm) node (v4)  {$v_4$}
		-- ++(180:3cm) node (v5)  {$v_5$}
		-- (v1);
		\draw (1,1.5) node (v6)  {$v_6$}
		-- (v1)
		-- ++(330:1.7cm) node (v7)  {$v_7$}
		-- ++(320:1.2cm) node (v8)  {$v_8$}
		-- ++(225:1.2cm) node (v9)  {$v_9$}
		-- (v5);
		\draw (4,1.5) node (v3)  {$v_3$}
		-- (v2)
		-- (v7)
		-- (v9)
		-- (v4)
		-- (v3);
		\draw (v6) -- (v5); }
	\node at (1.5,-1) {$H_1$};
	
	{\tikzstyle{every node}=[draw,circle,fill=white,minimum size=4pt,
		inner sep=0pt]
		\draw (6,3) node (v1)  {$v_1$}
		-- ++(0:3cm) node (v2)  {$v_2$}
		-- ++(270:3cm) node (v4)  {$v_4$}
		-- ++(180:3cm) node (v5)  {$v_5$};
		\draw (7,1.5) node (v6)  {$v_6$}
		-- (v1)
		-- ++(330:1.7cm) node (v7)  {$v_7$}
		-- ++(320:1.2cm) node (v8)  {$v_8$}
		-- ++(225:1.2cm) node (v9)  {$v_9$}
		-- (v5);
		\draw (10,1.5) node (v3)  {$v_3$}
		-- (v2)
		-- (v7)
		-- (v9)
		-- (v4)
		-- (v3);
		\draw (v6) -- (v5); }
	\node at (7.5,-1) {$H_2$};
	\node at (4.5,-2) {Fig. 3: $H_1$ and $H_2$};
	\end{tikzpicture}
\end{center}

Cerioli and Lima's work~\cite{CERIOLI2019} on showing that all longest paths have a common vertex in every $(K_{1,3}, P_5)$-free graph leads us to wonder if Gallai's question is true for graphs that forbid a Hamiltonian forbidden pair  $(R,S)$ as given in Theorem~\ref{forbiddenpairH}. The answer turns out to be positve and we obtain the following  result.

\begin{thm}\label{thm} Let   $R=K_{1,3}$, $S\in \{C_3,P_4,P_5,$ $P_6,Z_1,Z_2,Z_3,B_{1,1},B_{1,2}\}$, and $G$ be a connected $(R,S)$-free graph.
Then there exists a vertex common to all the  longest paths in $G$.
\end{thm}

We end this section by introducing some notation and terminologies.
Let $G$ be a graph. We use $V(G)$ and $E(G)$ to denote the vertex set and edge set of $G$, respectively. For two vertices $u,v\in V(G)$, we write $u\thicksim v$ if $uv\in E(G)$ and $u\nsim v$ otherwise. The set of neighbors of $u$ in $G$ is denoted by $N_G(u)$ or $N(u)$ if no confusion may arise. For any vertex $v\in V(G)$ and subset $X\subseteq V(G)$, the set of neighbors of $v$ in $X$ is denoted by $N_{X}(v)$ and $d_{X}(v):=|N_{X}(v)|$.
Let $X\subseteq V(G)$ be a vertex set, we use $G[X]$ to denote the subgraph induced by $X$ in $G$. A \emph{clique} in a graph $G$ is a subset of $V(G)$ that are pairwise adjacent.

A \emph{block} is a connected graph with no cutvertex, and a
\emph{block of $G$} is a maximal connected subgraph of $G$ that is
itself a block. Let $\mathcal{B}$ be the set of blocks and $\mathcal{C}$
the set of cutvertices of $G$. The {\it block-cutvertex tree}  of
a connected graph $G$
has vertex set $\mathcal{B}\cup \mathcal{C}$, and $c\in \mathcal{C} $
is adjacent to $B\in \mathcal{B}$  if and only if the block $B$
contains the cutvertex $c$.

A path  $P=v_{1}v_{2}\cdots v_{n}$ in graph $G$ is also called a $v_1v_{n}$-\emph{path}. For any two vertices $v_i,v_j\in V(P),i<j$, we use $v_iPv_j$ to denote the segment of $P$ starting at $v_i$ and ending at $v_j$, i.e., the subpath $v_{i}v_{i+1}\cdots v_{j}$. Let $P$ be a $uv$-path and $Q$ be an $xy$-path. If $P$ and $Q$ share a common vertex $w$, then we use $uPwQx$ to denote the concatenation of paths $uPw$ and $wQx$. For any two distinct vertices $u$ and $v$ in $G$, the \emph{distance} between  $u$ and $v$ in $G$, denoted by $d(u,v)$, is the length of a shortest $uv$-path in $G$. For any two disjoint subsets $X,Y\subseteq V(G)$, the \emph{distance} $d(X,Y):=\min\{d(u,v):u\in X,v\in Y\}$.

\section{Proof of Theorem~\ref{thm}}

In this section, we prove Theorem~\ref{thm}.   Our strategy is to first use the block-cutvertex tree of
$G$ to restrict the possible intersecting vertices for all longest paths within a block of $G$, and then apply the structural properties of the graph (properties implied by $(R,S)$-freeness) to find a common vertex to all longest paths in $G$.

\proof[Proof of Theorem~\ref{thm}]
Note that if $H_1\subseteq H_2$, then  $G$ is $H_1$-free implies that $G$ is also $H_2$-free.  Thus, to prove Theorem \ref{thm}, it suffices to consider the connected $(R,S)$-free graph $G$ when $(R,S)\in \{(K_{1,3}, P_6), (K_{1,3}, Z_3), (K_{1,3}, B_{1,2}), (K_{1,3}, N_{1,1,1})\}$.  By Theorem~\ref{forbiddenpairH1}, a connected $(K_{1,3}, N_{1,1,1})$-free graph is traceable, so every longest path contains all the vertices of $G$. Therefore, all longest paths have a common intersection
in a connected $(K_{1,3}, N_{1,1,1})$-free graph.  Hence, we only need to consider three pairs:
$(R,S)=(K_{1,3}, P_6), (K_{1,3}, Z_3)$, or $ (K_{1,3}, B_{1,2})$.
The proof when $(R,S)=(K_{1,3}, B_{1,2})$ is different from that for the other two pairs, so
we separate the remaining proof into two cases.

\medskip

{\noindent \bf Case 1: $(R,S)=(K_{1,3}, B_{1,2})$}.

\medskip


For this case, we use a known result by Furuya and Tsuchiya~\cite{FT2015} to reduce the problem to  a class of graphs called \emph{generalized comb}.
Thus, the problem regarding if all longest paths in a connected  $(K_{1,3}, B_{1,2})$-free graph have a common intersection lies in showing that Gallai's question is true for generalized combs. We now give the definition of a generalized comb, first introduced by Furuya and Tsuchiya~\cite{FT2015}.

A graph $G$ is a \emph{generalized comb} if it is obtained in the following way (see Fig. 4):
Let $m\geq3$ be an integer. Let $L_{i}(1\leq i\leq m)$ and $C$ be disjoint non-empty sets with
$|C|\geq m$, and let $R_i(1\leq i\leq m)$ be disjoint non-empty subsets of $C$. We define the
graph $G$ on  $\mathop{\cup}\limits_{1\leq i\leq m}L_{i}\cup C$ such that

$($i$)$ $L_{i}(1\leq i\leq m)$ and $C$ are cliques of $G$, and

$($ii$)$ for each $i(1\leq i\leq m)$, every vertex in $L_i$ is joined to all vertices in $R_i$.
Call $C$ the \emph{base} of the generalized comb $G$.

\begin{center}
\begin{tikzpicture}
\draw (0,0) ellipse (1 and 0.4);
\draw (3,0) ellipse (1 and 0.4);
\draw (6,0) ellipse (1 and 0.4);
\draw (0,2) ellipse (1 and 0.4);
\draw (3,2) ellipse (1 and 0.4);
\draw (6,2) ellipse (1 and 0.4);
\draw (-1,0) -- (-1,2);
\draw (1,0) -- (1,2);
\draw (2,0) -- (2,2);
\draw (4,0) -- (4,2);
\draw (5,0) -- (5,2);
\draw (7,0) -- (7,2);
\draw (3,0) ellipse (5 and 0.7) ;
\node at (0,2) {$L_1$};
\node at  (3,2) {$L_2$};
\node at  (6,2) {$L_m$};
\node at (0,0) {$R_1$};
\node at  (3,0) {$R_2$};
\node at  (6,0) {$R_m$};
\node at (0,1) {$+$};
\node at  (3,1) {$+$};
\node at  (6,1) {$+$};
\node at  (4.6,1) {\textbf{$\cdots$}};
\node at  (8,0.6) {$C$};
\node at (3,-1.5) {Fig. 4: Generalized Comb};
\end{tikzpicture}
\end{center}

Furuya and Tsuchiya~\cite{FT2015} showed that for a $(K_{1,3},B_{1,2})$-free graph $G$, the following hold.

$\mbox{{(i)}}$ If $G$ is connected, then $G$ is traceable unless $G$ is a generalized comb with at least three cutvertices.

$\mbox{{(ii)}}$ If $G$ is $2$-connected, then $G$ is Hamiltonian.

As Gallai's question is positive for traceable graphs, we can therefore assume that $G$ is a generalized comb with at least three cutvertices. By the construction of $G$, it is not hard to argue   that every longest path in $G$ contains all vertices in the base of the generalized comb $G$.
Thus, all longest paths in $G$ have a nonempty intersection.

\medskip

{\noindent \bf Case 2:  $(R,S)=(K_{1,3}, P_6)$ or $(R,S)=(K_{1,3}, Z_3)$}.

\medskip

 Let $\mathcal{B}$ be the set of blocks and $\mathcal{C}$
the set of cutvertices of $G$.
Let $T$ be the block-cutvertex tree of $G$ and $\mathcal{L}$ be the set of all longest paths in $G$.
 For each longest path $P\in \mathcal{L}$, define the set $$V_P=\{ B\in \mathcal{B}: V(P)\cap V(B)\ne \emptyset\}\cup \{c\in \mathcal{C} : c\in V(P) \},$$ and let
 $$T_P= T[V_P].$$
 Since $P$ is a connected subgraph of $G$, $T_P$ is a connected subgraph of $T$. Thus, $T_P$
 is a subtree of $T$.
 For any two longest paths $P$ and $Q$ in $\mathcal{L}$,  since $V(P)\cap V(Q)\ne \emptyset$,
 we have that
 $V(T_{P})\cap V(T_{Q})\neq \emptyset$. It is well known that a family of subtrees of a tree has Helly property (Let $\mathcal{F}$ be a family of sets. We say that $\mathcal{F}$ has \emph{Helly property} if and only if for every nonempty subfamily $\mathcal{H}\subseteq \mathcal{F}$ and for all sets $X,Y\in \mathcal{H}$ such that $X\cap Y\neq \emptyset$, then $\mathop{\cap}\limits_{X\in \mathcal{H}}X\neq \emptyset$. See problem 18 on page 49 of \cite{Lovasz1993}), so there is a vertex $B\in V(T)$ such that $B\in \mathop{\cap}\limits_{P\in \mathcal{L}}V(T_{P})$. By the construction of $T$, $B$ is either a cutvertex or a block of $G$. If $B$ is a cutvertex of $G$, then $B$ is a common vertex of all longest paths in $G$.  Thus, we assume that $B$ is a block. So $B=K_2$ is an edge or is a 2-connected graph.  If $B=xy$ is an edge, then $B$ is a cutedge of $G$. We claim that either $x$ or $y$ is a common vertex of all longest paths in $G$. Otherwise, there would be two longest paths $P,Q\in \mathcal{L}$ such that $x\in V(P)$ but $y\not\in V(P)$ and $y\in V(Q)$ but $x\not\in V(Q)$. Since $V(P)\cap V(Q)\neq \emptyset$, there exists an $xy$-path in $P\cup Q$. This contradicts  the fact that $xy$ is a cutedge of $G$. Therefore, we assume that $B$ is a 2-connected subgraph of $G$.
 By Theorems~\ref{forbiddenpairH1} and~\ref{CZ3}, $B$ has a Hamiltonian cycle or Hamiltonian path.  (The two graphs $H_1$ and $H_2$ in Theorem~\ref{CZ3} have a Hamiltonian path by inspections.)

Let $\mathcal{P}$ be the set of all longest paths of $G$ that does not contain all vertices of $B$.
We may assume that  $\mathcal{P} \ne \emptyset$. For otherwise, we are done by noting that  $B$ has a Hamiltonian path and every longest path contains all vertices of $B$.   This assumption that  $\mathcal{P} \ne \emptyset$ implies also that  $B$ contains
a cutvertex of $G$.
For a path $P\in \mathcal{P}$, let $u_p$ be an endvertex of $P$ that is not contained in $B$, and let
    $x_p$ be  the vertex of $P$ that is closest to $u_p$ on $P$. We call $x_pPu_p$
    a \emph{pendent segment} of $P$ on $B$. Clearly,  $x_p$ is a cutvertex of $G$ and $V(x_pPu_p)\cap V(B)=\{x_p\}$.  The vertex $x_p$ is called an \emph{attachment} of $P$ on $B$.

\begin{claim}\label{claim0}  The following statements hold.
	\begin{enumerate}[$(i)$]
		\item For every longest path $P$ in $G$, $|V(P)|\ge |V(B)|$; furthermore, $|V(P)|\ge |V(B)|+1$
		if $B$ is Hamiltonian.
		\item For every $P\in \mathcal{P_{}}$, $P$ has exactly  two different pendent segments on $B$.
	\end{enumerate}

\end{claim}

\proof[Proof of Claim \ref{claim0}] Since $B$ has a Hamiltonian path, every longest path of $G$
contains at least $|V(B)|$ vertices. When $B$ is Hamiltonian,  let $C$ be a Hamiltonian cycle of $B$. For any $z\in V(B)$ such that $z$ is a cutvertex of $G$, let $z_1$ be a neighbor of $z$ in $G$
from $V(G)\setminus V(B)$.  Let $z^*$ be a neighbor of $z$ on $C$. Then $z^*Czz_1$ is a path of $G$
that contains $|V(B)|+1$ vertices. This proves (i).

For (ii), note that every attachment  of $P$ on $B$ is a cutvertex of $G$. If $P$
has exactly one pendent segment on $B$, then by the second part of the proof for (i) above, we know that $V(B)\subseteq V(P)$, contradicting to the definition of $\mathcal{P}$.
Since $P$ is a path and its attachment on $B$ is a cutvertex of $G$,
$P$ has exactly  two different pendent segments on $B$.
\qed

For every $P\in \mathcal{P_{}}$, we use $x_{p}$ to denote an attachment of $P$ on  $B$, and let $y_{p}$ be a neighbor of $x_p$ in $G$ from $V(G)\setminus V(B)$.

\begin{claim}\label{claim1}  The following statements hold.
	\begin{enumerate}[$(i)$]
	
	\item  For any vertex $x\in V(B)$, $d_{B}(x)\geq 2;$

\item  $N_{B}(x_{p})$ is a clique$;$ and

\item $N_{B}(x_{p})\subseteq V(P)$.
	\end{enumerate}

\end{claim}

\proof[Proof of Claim~\ref{claim1}](i) is clear  since $B$ is 2-connected.
(ii) follows from the fact that   $G$ is $K_{1,3}$-free and $x_p$ is a cutvertex of $G$.
For (iii),
suppose that there exists a vertex $z_{p}\in N_{B}(x_{p})\setminus V(P)$. Let $w_{p}\in V(B)$ be the immediate successor of vertex $x_{p}$ on $P$. By (ii), $z_{p}\thicksim w_{p}$, then replacing the edge $x_{p}w_{p}$ by $x_{p}z_{p}w_{p}$ gives  a path in $G$ that is longer than $P$. This gives  a contradiction to the assumption  that $P$ is a longest path in $G$.
\qed

We now separate the proof of Case 2 into two subcases.

\medskip

{\bf \noindent Subcase 2.1 $(R,S)=(K_{1,3}, P_6)$}.

\medskip

\begin{claim}\label{claim3} Let $P,Q\in \mathcal{P}$. Then $N_{B}(x_{p})\cap N_{B}(x_{q})\neq \emptyset$.
\end{claim}

\proof [Proof of Claim~\ref{claim3}] If $x_{p}\thicksim x_{q}$ or $x_{p}=x_{q}$, then Claim \ref{claim1} implies that $N_{B}(x_{p})\cap N_{B}(x_{q})\neq \emptyset$. So we assume that $x_{p}\nsim x_{q}$ and $x_{p}\ne x_{q}$. Suppose that $N_{B}(x_{p})\cap N_{B}(x_{q})=\emptyset$. Then $d(x_{p},x_{q})\geq 3$. Let $Q'$ be a shortest $x_{p}x_{q}$-path in $B$. Then $y_px_pQ'x_qy_q$ contains  an induced $P_6$ in $G$. This gives a contradiction to the assumption that $G$ is $P_6$-free. (See Fig. 5 for an illustration of the $P_6$.)
		\qed

\begin{center}
\begin{tikzpicture}
\draw (3,0) ellipse (2 and 1);

{\tikzstyle{every node}=[draw,circle,fill=black,minimum size=4pt,
                            inner sep=0pt]
\draw (0.3,-1.5)
        --(0.8,-1) node (y1)  { }
        -- (1.3,-0.5) node (x1)  { }
        -- (4.7,-0.5) node (x2)  { }
        -- (5.2,-1) node (y2)  { }
        -- (5.7,-1.5)   { };
 }
\node at (0.4,-1) {$y_{p}$};
\node at (0.9,-0.5) {$x_{p}$};
\node at (5.1,-0.5) {$x_{q}$};
\node at (5.6,-1) {$y_{q}$};
\node at (3,-0.2){$Q'$};
\node at (4.8,0.8) {$B$};
\node at (3,-2) {Fig. 5: Illustration of Claim \ref{claim3}};
\end{tikzpicture}
\end{center}

\begin{claim}\label{claim4} $\mathop{\cap}\limits_{P\in\mathcal{P}_{}}N_{B}(x_{p})\neq \emptyset$.
\end{claim}

\proof[Proof of Claim~\ref{claim4}] We prove by contradiction.
Let $\ell$ be the smallest index such that there exists
$\ell$ paths, say $P_1, P_2, \cdots, P_\ell$ from $\mathcal{P}$ so that
\begin{equation}\label{nocommon}
\mathop{\cap}\limits_{1\le  i\le \ell}N_{B}(x_{p_i})=\emptyset.
\end{equation}
By Claim~\ref{claim3}, $\ell\ge 3$.  By the choice of $\ell$,
$$\mathop{\cap}\limits_{1\le  i\le \ell-1}N_{B}(x_{p_i})\ne \emptyset,\quad \mbox{and}\quad\mathop{\cap}\limits_{2\le  i\le \ell}N_{B}(x_{p_i})\ne \emptyset.$$
 Let $$z_1\in \mathop{\cap}\limits_{1\le  i\le \ell-1}N_{B}(x_{p_i}), \quad \mbox{and}\quad  z_2\in \mathop{\cap}\limits_{2\le  i\le \ell}N_{B}(x_{p_i}).$$
 By the assumption in~\eqref{nocommon}, we have that $z_1\ne z_2$.

 We claim that $x_{p_1}\nsim z_2$, $x_{p_\ell}\nsim z_1$, and $x_{p_1}\nsim x_{p_\ell}$.
 If $x_{p_1}\sim z_2$, i.e., $z_2\in N_B(x_{p_1})$, then $z_2\in \mathop{\cap}\limits_{1\le  i\le \ell}N_{B}(x_{p_i})$, showing a contradiction to the choice of $\ell$. Similarly, $x_{p_\ell}\nsim z_1$.  If $x_{p_1}\sim x_{p_\ell}$, then $x_{p_\ell}\sim z_1$.
 This is because $x_{p_1}\sim z_1$ and $N_B(x_{p_1})$ is a clique by  Claim~\ref{claim1} (ii).

 Since $z_1$ and $z_2$ are both adjacent to $x_{p_2}$, we have that $z_1\thicksim z_2$   by Claim~\ref{claim1} (ii). Then $y_{p_1}x_{p_1}z_1z_2x_{p_{\ell}}y_{p_{\ell}}$ is an induced $P_6$ in $G$.
 This gives a  contradiction to the $P_6$-freeness assumption of $G$. (See Fig. 6 for an illustration.)
 \qed
\begin{center}
\begin{tikzpicture}
\draw (3,0) ellipse (2 and 1);

{\tikzstyle{every node}=[draw,circle,fill=black,minimum size=4pt,
                            inner sep=0pt]
\draw (0.3,-1.5)
        --(0.8,-1) node (y1)  { }
        -- (1.3,-0.5) node (x1)  { }
        --(1.7,0.8) node (z1)  { }
        --(2.2,-0.9) node (x2)  { }
        --(3.5,1) node (z2)  { }
        --(3.1,-1) node (xi)  { }
        -- (z2)   { }
        --(4.4,-0.7) node (xi1)  { }
        --(4.8,-1.2) node (yi1)  { }
        --(5.4,-1.9)   { };
 }

\draw (xi)--(z1)--(z2);
\node at (0.4,-1) {$y_{p_1}$};
\node at (0.9,-0.5) {$x_{p_1}$};
\node at (1.7,1.1) {$z_1$};
\node at (2,-1.2) {$x_{p_2}$};
\node at (3.5,1.2) {$z_2$};
\node at (3.3,-1.3) {$x_{p_{\ell-1}}$};
\node at (4.9,-0.8) {$x_{p_{\ell}}$};
\node at (5.3,-1.3) {$y_{p_{\ell}}$};
\node at (2.6,-1.2){$\cdots$};
\node at (4.8,0.8) {$B$};
\node at (3,-2.5) {Fig. 6: Illustration of Claim \ref{claim4}};
\end{tikzpicture}
\end{center}

Since any longest path $P$ with $P\not\in \mathcal{P}_{}$ contains all vertices of $B$,  Claim \ref{claim1} (iii) and Claim \ref{claim4} imply that all longest paths in $G$ have a nonempty intersection.

\medskip

{\bf \noindent Subcase 2.2 $(R,S)=(K_{1,3}, Z_3)$}.

\medskip

If $B$ is isomorphic to $H_1$, by Claim \ref{claim1} (ii), we know that the attachments of any longest path $P\in \mathcal{P}$ can only be chosen from $\{v_3,v_6,v_8\}$.  For any $i\ne j$, $i,j\in \{3,6,8\}$,
every longest $v_iv_j$-path in $B$ contains all vertices of $B$. By Claim~\ref{claim0},
any $P\in \mathcal{P}$ has exactly two attachments on $B$.
It then follows that $V(B)\subseteq V(P)$ for any longest path $P\in\mathcal{P}$, contradicting to the choice of $\mathcal{P}$. Therefore, $B$ is not isomorphic to $H_1$. Similarly, $B$ is not isomorphic to $H_2$. Therefore, $B$ has at least 10 vertices and thus is Hamiltonian by Theorem~\ref{forbiddenpairH}.

\begin{claim}\label{3vertices-out-B} We may assume that for each longest path $P\in \mathcal{P}$, $P$ contains at least three vertices outside of $B$, i. e., $|V(P)\cap(V(G)\setminus V(B))|\geq 3$.
\end{claim}

\proof[Proof of Claim~\ref{3vertices-out-B}]Recall that for every $P\in \mathcal{P}$, $P$ contains two different segments which do not belong to $B$. Therefore, $|V(P)\cap(V(G)\setminus V(B))|\geq 2$. Suppose that there exists $P\in \mathcal{P}$ such that $|V(P)\cap(V(G)\setminus V(B))|=2$. Then $|V(P)|\leq |V(B)|+1$, since $P$ does not contain all vertices of $B$ by the definition of $\mathcal{P}$.
By Claim~\ref{claim0} (i), we then know that $|V(P)|=|V(B)|+1$.
 Then $P$ contains all except one vertex of $B$.

 Let $P\in \mathcal{P}$ be a fixed path with an attachment $x_p$. Let
 $$\mathcal{P}^{1}=\{Q\in \mathcal{P}: Q\ \mbox{contains}\ x_{p}\}, \quad \mbox{and}\quad \mathcal{P}^{2}=\{Q\in \mathcal{P}: Q\ \mbox{does\ not\ contain}\ x_{p}\}.$$

Note that  $N_{B}(x_{p})$ is a clique of size at least 2 in $B$  by Claim~\ref{claim1} (i) and (ii).
Thus, for each path $Q\in \mathcal{P}^{1}$, the fact that $Q$ contains $x_p$ and misses exactly one vertex of $B$
implies that $N_{B}(x_{p})\subseteq V(Q)$. For any path $Q\in \mathcal{P}^{2}$, $N_{B}(x_{p})\subseteq V(Q)$ since $Q$ contains all except one vertex of $B$. Therefore, $N_{B}(x_{p})\subseteq \mathop{\cap}\limits_{P\in\mathcal{P}}V(P)$.
  This together with the fact that every longest path $P$ of $G$ with $P\not\in \mathcal{P}_{}$ contains all vertices of $B$, we see that all longest paths in $G$ have a nonempty intersection.
     Therefore, we may assume that $|V(P)\cap(V(G)\setminus V(B))|\geq 3$.
\qed

By Claim~\ref{3vertices-out-B}, for each  path $P\in \mathcal{P}$,  $|V(P)\cap(V(G)\setminus V(B))|\geq 3$. Without loss of generality, we assume that $x_{p}$ is an attachment of $P$ such that the segment of $P$ attached with $x_{p}$ contains at least two vertices outside of $B$ when $P$ is considered in the following.

\begin{claim}\label{claim5} Let $P,Q\in \mathcal{P}$. Then $N_{B}(x_{p})\cap N_{B}(x_{q})\neq \emptyset$.
\end{claim}

\proof[Proof of Claim~\ref{claim5}]If $x_{p}\thicksim x_{q}$ or $x_{p}=x_{q}$, then Claim \ref{claim1} implies that $N_{B}(x_{p})\cap N_{B}(x_{q})\neq \emptyset$. So we assume that $x_{p}\ne x_{q}$ and
	$x_{p}\nsim x_{q}$. Suppose that $N_{B}(x_{p})\cap N_{B}(x_{q})=\emptyset$. Let $w_{p}\in N_{B}(x_{p})$  and $w_{q}\in N_{B}(x_{q})$
		such that $d(w_{p},w_{q})=d(N_{B}(x_{p}),N_{B}(x_{q}))\geq 1$.
		Let $z_p\in N_{B}(x_{p})$ with $z_p\ne w_p$, and $z_q\in N_{B}(x_{q})$ with $z_q\ne w_q$.
		By the assumption that  $N_{B}(x_{p})\cap N_{B}(x_{q})=\emptyset$,
		we know that
		$$
		x_p\not\sim w_q, x_q\not\sim w_p.
		$$

		Take $Q'$ be a shortest $w_{p}w_{q}$-path in $B$.
								We consider two subcases below according to the
		length of $Q'$.

\medskip
{\bf \noindent Subcase $|V(Q')|\geq 3$.}

\medskip
Let $u_p$ be the neighbor of $w_p$ on $Q'$.
Note that $z_{p}$ is adjacent to $w_p$ and possibly $u_p$ from $Q'$ but nothing else
by the  choice of $w_{p}$.
If $z_p\sim u_p$, then  $G[\{z_p, x_q,y_q\}\cup V(Q')]$
contains an induced $Z_3$ in $G$., where $z_pw_pu_p$ forms the triangle of the $Z_3$.
If $z_p\not\sim u_p$, then $G[\{z_{p},x_{p},x_{q}\}\cup V(Q')]$ contains an induced  $Z_3$ in $G$,
where $z_px_pw_p$ forms the triangle.
These give a contradiction to the $Z_3$-freeness assumption of $G$. (See Fig. 7 for an illustration.)

\begin{center}
\begin{tikzpicture}
\draw (3,0) ellipse (2 and 1);

{\tikzstyle{every node}=[draw,circle,fill=black,minimum size=4pt,
                            inner sep=0pt]
\draw (0.3,-1.5)
        --(0.8,-1) node (y1)  { }
        -- (1.3,-0.5) node (x1)  { }
        --(2.2,-0.9) node (z1)  { }
        --(1.7,0.8) node (w1)  { }
        --(3.5,1) node (w2)  { }
        --(3.1,-1) node (z2)  { }
        --(4.4,-0.7) node (x2)  { }
        --(4.8,-1.2) node (y2)  { }
        --(5.4,-1.9)   { };
 }
\draw (x1)--(w1);
\draw (x2)--(w2);
\node at (0.4,-1) {$y_{p}$};
\node at (0.9,-0.5) {$x_{p}$};
\node at (1.7,1.1) {$w_{p}$};
\node at (2,-1.2) {$z_{p}$};
\node at (3.5,1.2) {$w_{q}$};
\node at (3.3,-1.3) {$z_{q}$};
\node at (4.9,-0.8) {$x_{q}$};
\node at (5.3,-1.3) {$y_{q}$};
\node at (2.7,0.6){$Q'$};
\node at (4.8,0.8) {$B$};
\node at (3,-2.5) {Fig. 7: Illustration of Claim \ref{claim5}};
\end{tikzpicture}
\end{center}

\medskip

{\bf \noindent Subcase  $|V(Q')|=2$, i.e., $Q'=w_{p}w_{q}$.}

\medskip

If $x_{p}$ has at least two neighbors in $V(G)\setminus V(B)$, say $y_{p},v_{p}$, then $G[\{x_{p},y_{p},v_{p}\}]$ is a triangle (otherwise, we find an induced $K_{1,3}$ in $G$). Since $ x_{p}\nsim x_{q}, , x_{p}\nsim w_{q},x_{q}\nsim w_{p}$, $G[\{x_{p},y_{p},v_{p},w_{p},w_{q},x_{q}\}]$ is an induced $Z_3$ in $G$, a contradiction. Therefore, $x_{p}$ has exactly one neighbor $y_{p}$ in $V(G)\setminus V(B)$. By the assumption on $x_p$, $y_{p}$ has a neighbor $v_{p}\in V(G)\setminus V(B)$ and $x_p\nsim v_p$. If $z_{q}\nsim w_{p}$, then $G[\{x_{q},z_{q},w_{q},w_{p},x_{p},y_{p}\}]$ is an induced $Z_3$ in $G$, a contradiction. So  $z_{q}\thicksim w_{p}$, then $G[\{z_{q},w_{p},w_{q},x_{p},y_{p},v_{p}\}]$ is an induced  $Z_3$ in $G$, a contradiction.
\qed

\begin{claim}\label{claim6} $\mathop{\cap}\limits_{P\in\mathcal{P}}N_{B}(x_{p})\neq \emptyset$.
\end{claim}

\proof[Proof of Claim~\ref{claim6}]
We prove by contradiction.
Let $\ell$ be the smallest index such that there exist
$\ell$ paths, say $P_1, P_2, \cdots, P_\ell$ from $\mathcal{P}$ so that
\begin{equation}\label{nocommon2}
\mathop{\cap}\limits_{1\le  i\le \ell}N_{B}(x_{p_i})=\emptyset.
\end{equation}
By Claim~\ref{claim3}, $\ell\ge 3$.  By the choice of $\ell$,
$$\mathop{\cap}\limits_{1\le  i\le \ell-1}N_{B}(x_{p_i})\ne \emptyset,\quad \mbox{and}\quad\mathop{\cap}\limits_{2\le  i\le \ell}N_{B}(x_{p_i})\ne \emptyset.$$
Let $$z_1\in \mathop{\cap}\limits_{1\le  i\le \ell-1}N_{B}(x_{p_i}), \quad \mbox{and}\quad  z_2\in \mathop{\cap}\limits_{2\le  i\le \ell}N_{B}(x_{p_i}).$$
By the assumption in~\eqref{nocommon2}, we have that $z_1\ne z_2$.
Since $z_1$ and $z_2$ are both adjacent to $x_{p_2}$, we have that $z_1\thicksim z_2$   by Claim~\ref{claim1} (ii).

If $x_{p_1}$ has at least two neighbors in $V(G)\setminus V(B)$, say $y_{p_1},u_{p_1}$, then $G[\{x_{p_1},y_{p_1},u_{p_1}\}]$ is a triangle. Since $x_{p_1}\nsim z_{2}, x_{p_{\ell}}\nsim z_1, x_{p_1}\nsim x_{p_{\ell}}$, $G[\{x_{p_1},y_{p_1},u_{p_1},z_{1},z_{2},x_{p_{\ell}}]$ is an induced $Z_3$ in $G$, a contradiction. Therefore, $x_{p_1}$ has exactly one neighbor $y_{p_1}$ in $V(G)\setminus V(B)$. By the assumption on $x_p$, $y_{p_1}$ has a neighbor $u_{p_1} \in V(G)\setminus V(B)$ and $x_{p_1}\nsim u_{p_1}$. Then $G[\{x_{p_2},z_1,z_2,x_{p_1},y_{p_{1}},u_{p_{1}}\}]$ is an induced $Z_3$ in $G$, a contradiction.
\qed

Since every longest path $P$ with $P\not\in \mathcal{P}$ contains all vertices of $B$,  Claim \ref{claim1} (iii) and Claim \ref{claim6} imply that all longest paths in $G$ have a nonempty intersection.
This finishes the proof of Subcase 2.2.
\qed

\section{Conclusion and open problems}

In this paper, we showed that Gallai's question is true for connected $(R,S)$-free graphs, where  $R=K_{1,3}$, and $S\in \{C_3,P_4,P_5,P_6,Z_1,Z_2,Z_3,B_{1,1},B_{1,2}\}$. It is natural to consider whether Gallai's question is positive for the superclasses of the above graphs, such as $K_{1,3}$-free graphs, $C_3$-free graphs, $P_6$-free graphs, etc..
 It is still unknown whether every set of three longest paths of a connected graph share a common vertex, and this is  conjectured to be true in \cite{HHM2008}.
 We suspect the following questions might have a positive answer.

 \begin{prob}
 	Do all longest paths have a nonempty intersection in a connected $K_{1,3}$-free or $P_5$-free graph?
 \end{prob}

 \begin{prob}
	Does every set of three longest paths have a nonempty intersection in a connected $K_{1,3}$-free or $P_5$-free graph?
\end{prob}		

\bibliography{Forbidden-pairs-for-nonempty-intersection-of-longest-paths}

\end{document}